\documentclass[12pt,reqno]{amsart}
\usepackage{amsmath}
\usepackage{amssymb}
\usepackage{amsthm}
\usepackage[left=2.9cm,top=2.5cm,right=2.9cm,bottom=2.6cm]{geometry}
\usepackage{epsfig}
\usepackage[cp1250]{inputenc}
\usepackage[T1]{fontenc}
\usepackage{txfonts}
\usepackage{verbatim}

\begin{document}
\newtheorem{theorem}{Theorem}[section]
\newtheorem{lemma}[theorem]{Lemma}
\newtheorem{definition}[theorem]{Definition}
\newtheorem{conjecture}[theorem]{Conjecture}
\newtheorem{proposition}[theorem]{Proposition}
\newtheorem{algorithm}[theorem]{Algorithm}
\newtheorem{corollary}[theorem]{Corollary}
\newtheorem{observation}[theorem]{Observation}
\newtheorem{problem}[theorem]{Problem}
\newtheorem{example}[theorem]{Example}
\newcommand{\noin}{\noindent}
\newcommand{\ind}{\indent}
\newcommand{\al}{\alpha}
\newcommand{\om}{\omega}
\newcommand{\pp}{\mathcal P}
\newcommand{\ppp}{\mathfrak P}
\newcommand{\R}{{\mathbb R}}
\newcommand{\N}{{\mathbb N}}
\newcommand{\Z}{{\mathbb Z}}
\newcommand{\ON}{\mathrm{ON}}
\newcommand\eps{\varepsilon}
\newcommand{\E}{\mathbb E}
\newcommand{\Prob}{\mathbb{P}}
\newcommand{\pl}{\textrm{C}}
\newcommand{\rk}{\textrm{rk}}
\newcommand{\cro}{\kappaup}
\newcommand{\dist}{\textrm{dist}}
\newcommand{\rad}{\textrm{rad}}
\newcommand{\diam}{\textrm{diam}}
\renewcommand{\labelenumi}{(\roman{enumi})}
\newcommand{\bc}{\bar c}
\newcommand{\G}{{\mathfrak S}}
\newcommand{\T}{{\mathfrak T}}
\newcommand{\Remark}[1]{\par \noindent \textnormal{\textbf{Note.~}} #1}

\title{Cops and Robbers ordinals of cop-win trees}
\author{Anthony Bonato}
\address{Ryerson University}
\email{\tt abonato@ryerson.ca}
\author{Przemys\l{}aw Gordinowicz}
\address{Institute of Mathematics, Lodz University of Technology, \L{}\'od\'z, Poland}
\email{\tt pgordin@p.lodz.pl}
\author{Ge\v na Hahn}
\address{Université de Montréal}
\email{\tt hahn@iro.umontreal.ca}

\thanks{The first and third authors gratefully acknowledge support from NSERC}
\subjclass{05C63,05C57}
\begin{abstract}
A relational characterization of cop-win graphs was provided by
Nowakowski and Winkler in their seminal paper on the game of Cops and
Robbers. As a by-product of that characterization, each cop-win graph
is assigned a unique ordinal, which we refer to as a CR-ordinal. For
finite graphs, CR-ordinals correspond to the length of the game
assuming optimal play, with the cop beginning the game in a least
favourable initial position. For infinite graphs, however, the
possible values of CR-ordinals have not been considered in the
literature until the present work.

We classify the CR-ordinals of cop-win trees as either a finite
ordinal, or those of the form $\alpha + \omega$, where $\alpha$ is a
limit ordinal. For general infinite cop-win graphs, we provide an
example whose CR-ordinal is not of this form. We finish with some
problems on characterizing the CR-ordinals in the general case of
cop-win graphs.
\end{abstract}
\keywords{Cops and Robbers, cop-win graphs, trees, ordinal numbers}


\vspace{150pt}

\maketitle

\section{Introduction\label{sec:intro}}
In the game of Cops and Robbers, cop-win graphs are those for which one cop has a winning strategy to capture the robber.  Nowakowski and Winkler~\cite{nw} were the first ones to provide a
characterization of cop-win graphs that is not restricted to finite graphs.  Consider an infinite graph $G$ of order $\aleph_\gamma$ and write $\kappa=\aleph_{\gamma+1}$, where $\gamma \ge 0$ is an
ordinal. For a vertex $v$, let $N[v]$ be its closed neighbourhood. Define the relations $\{\le_\alpha\}_{\alpha < \kappa}$ on $V(G)$ as follows (it is useful to remember that  $\leq_\alpha\subseteq
V(G)\times V(G)$).
\begin{enumerate}
\item If $u=v$, then $u \le_0 v$.
\item $u \le_\alpha v \textrm{~if
for all~} x \in N[u], \textrm{~there exists~} y \in N[v] \textrm{~such
that~} x \le_\beta y \textrm{~for some~} \beta < \alpha$.
\end{enumerate}

Observe that for all $\alpha\ <\ \beta$, $\le_\alpha\ \subseteq\ \le_\beta$. It follows that this tower of relations stabilizes at some ordinal $\rho\leq \kappa$; that is, there is a minimum $\rho$
such that $\le_\rho\ =\ \le_{\rho+1}$. We refer to $\rho$ as the \emph{CR-ordinal} of $G$. In the case $G$ is a finite graph of order $n$, note that $\rho \leq n(n-1).$

A result of \cite{nw} that had often been overlooked is that the graph $G$ is cop-win if and only if the relation $\le_\kappa$ is trivial (that is, $\le_{\kappa} = V(G)\times V(G))$. In this paper,
we refer to $\le_\rho$ as the \emph{capture relation} on cop-win graphs. The capture relation has provided insights into algorithms for recognizing cop-win graphs \cite{hm}, graphs with higher cop
number \cite{cm}, and for generalized Cops and Robbers games \cite{bm}.  It is not well-defined for graphs whose cop number is greater than one and that is one reason we still do not have a good
characterization of $k$-cop-win graphs, where $k>0.$

We say that an ordinal $\kappa$ is a \emph{CR-ordinal} if it is the CR-ordinal of some cop-win graph. We denote the CR-ordinal of a given graph $G$ by $\rho(G)$. A basic question, therefore, is which
ordinals are CR-ordinals?

For finite graphs, $\rho(G)$ is an integer and is the
length of the game assuming that the cop and robber play optimally
(that is, the cop plays to minimize the length of the game, while the
robber plays to maximize it), \emph{maximized} over all possible
starting positions for the cop. A related notion is the \emph{capture
time} of $G$, which is the length of the game assuming optimal play
\emph{minimized} over starting positions for the cop; see
\cite{boncapt}. In the infinite case, we can no longer directly
interpret the CR-ordinal as a length of a game. We emphasize that even
in the infinite case, a cop captures a robber in a cop-win graph in
finitely many rounds as there is no infinite descending
sequence of ordinals.

Finite paths demonstrate that every non-negative integer can be a CR-ordinal. Not surprisingly, the classification of CR-ordinals for infinite cop-win graphs is more complex. Indeed, \cite{bht}
suggests that the family of cop-win graphs is likely not classifiable as for each infinite cardinal $\aleph_\gamma$ there are the maximum possible number $2^{\aleph_\gamma}$ of non-isomorphic cop-win
graphs of order $\aleph_\gamma$.

We provide a characterization of CR-ordinals for cop-win trees in Theorem~\ref{thm:tree}. The cop-win trees are precisely the \emph{rayless ones}: trees not containing an infinite path as a subgraph.
Theorem~\ref{thm:tree} shows that finite ordinals and the ordinals of the form $\alpha+\omega$, where $\alpha$ is a limit ordinal and $\omega$ is the first infinite ordinal, are the CR-ordinals
witnessed by cop-win trees. In the final section, we consider general cop-win graphs. We provide a family of cop-win graphs, inspired by the graphs in \cite{polat}, whose CR-ordinals are not those
found from trees. Some open problems on the classification of CR-ordinals are stated at the end of the paper.

All the graphs we consider are simple and undirected.   We assume the reader is (or will be) familiar with the basic game of Cops and Robbers as defined in, for example,~\cite{bonato,nw}. For
additional background on Cops and Robbers and its variants, see the book~\cite{bonato} and the surveys~\cite{bonato1,bonato2,bonato3,hahn}. For background on graph theory, see \cite{diestel,west}. We
denote the distance between vertices $u$ and $v$ by $d(u,v)$. If $u$ is a vertex and $S$ a set of vertices, then $d(u,S)$ is the minimum distance from $u$ to a vertex in $S$. Let $\mathrm{ON}$ be the
proper class of ordinals. We use the property that every ordinal is the set of ordinals preceding it in the well ordering of $\ON$. For instance, $\omega$ consists of the set of all finite ordinals,
and we use this notation throughout. Hence, $\omega =\mathbb N =\{0,1,2,\ldots\}$. Recall that a \emph{successor ordinal} is one which contains a maximum element; such ordinals are of the form
$\alpha+1$, where $\alpha$ is some ordinal. A \emph{limit ordinal} is not a successor ordinal; for example, $\omega$ is a limit ordinal. Transfinite induction is analogous to usual induction, but
considers the cases of both successor and limit ordinals.  For further reading on ordinals and cardinals, see \cite{devlin,ros}.

\section{Capture-time ordinal}

Throughout this section, let $G$ be a cop-win graph. Before we state
our main result in the next section, it will be useful to use the
sequence $\{\leq_\alpha\}_{\alpha\leq \rho(G)}$ to introduce a
parameter  that provides a simpler means of
computing $\rho(G)$. Capture-time \cite{boncapt} is a temporal
counterpart to the cop number for a graph, measuring the length of the
game assuming optimal play. We now provide an ordinal analogue of
capture-time.

\begin{figure}[htbp]
\begin{minipage}{.4\textwidth}
\begin{center}
\epsfig{clip, figure=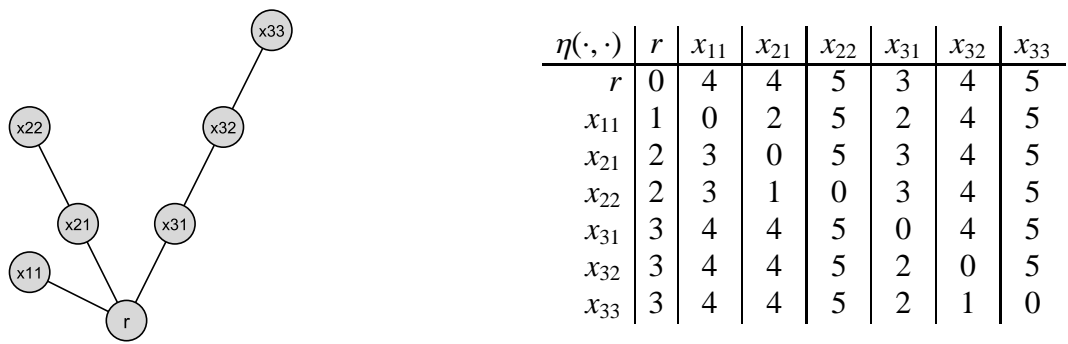, width=.7\textwidth}
\end{center}
\end{minipage}
\hfill
\begin{minipage}{.5\textwidth}
\begin{center}
\begin{tabular}{r|*{6}{c|}c}
$\eta(\cdot, \cdot)$ & $r$ & $x_{11}$ & $x_{21}$ & $x_{22}$ &
$x_{31}$ & $x_{32}$ & $x_{33}$ \\
\hline
$r$ & 0 & 4 & 4 & 5 & 3 & 4 & 5 \\
$x_{11}$ & 1 & 0 & 2 & 5 & 2 & 4 & 5 \\
$x_{21}$ & 2 & 3 & 0 & 5 & 3 & 4 & 5 \\
$x_{22}$ & 2 & 3 & 1 & 0 & 3 & 4 & 5 \\
$x_{31}$ & 3 & 4 & 4 & 5 & 0 & 4 & 5 \\
$x_{32}$ & 3 & 4 & 4 & 5 & 2 & 0 & 5 \\
$x_{33}$ & 3 & 4 & 4 & 5 & 2 & 1 & 0 \\
\end{tabular}
\end{center}
\end{minipage}
\caption{The tree $T$ and its table of $\eta$ values.}\label{fig:FiniteTree}
\end{figure}

For $u,v \in V(G)$, define $\eta(u,v) = \alpha$, where $\alpha$ is the minimum ordinal for which $u \le_\alpha v$ holds. Note that $\eta(u,v)$ is well-defined as ordinals are well-ordered. If
$\eta(u,v)$ is finite, then we may interpret it as the length of time it takes a cop on $v$ to capture a robber on $u$, assuming both play optimally and the robber moves first. Note that the relation
is not necessarily symmetric: $\eta(u,v)$ may be different than $\eta(v, u)$. For an example, see Figure~\ref{fig:FiniteTree} and its corresponding table of $\eta$ values.

Define $\eta(v)$ as the minimum ordinal $\alpha$ such that $u
\le_\alpha v$ holds for every $u \in V(G)$. Observe that from the
definitions, for any $v \in V(G)$ we have that $\eta(v) = \sup_{u \in
V(G)} \eta(u, v)$. Finally, we define $$\eta(G) = \min_{v \in V(G)}
\eta(v).$$ When finite, $\eta(G)$ is precisely the capture time
\cite{boncapt} of the cop-win graph $G$, and, hence, we will call such
ordinals the \emph{capture-time ordinals} associated with cop-win
graphs.

A crucial observation (which follows from the definitions) is that
$$\rho(G) =\sup_{v \in V(G)} \eta(v).$$ In particular, in the finite
case, $\rho(G)$ is the maximum capture time over all initial positions
of the cop.

Define $\theta(G)$ to be the set of vertices that realize $\eta(G);$
namely $$\theta(G) = \{v \in V(G) \colon \eta(v) = \eta(G)\}.$$ Note
that by the definitions, $\theta (G) \neq \emptyset.$ We may view the
set $\theta(G)$ as the set of vertices which are optimal starting positions
for the cop. For example, in a finite tree $T$, $\theta(T)$ is the
centre of the tree.

As we have just introduced a number of graph parameters, we give an
example that illustrates them. Consider the tree $T$
depicted in Figure~\ref{fig:FiniteTree}, along with its table of
$\eta$-values.
Note that by considering the table, we derive that $\rho(T) = 5$, $\eta(T)=3$, and $\theta(T) =
\{ r, x_{31} \}$.

\smallskip

\begin{figure}[htbp]
\begin{center}
\epsfig{clip, figure=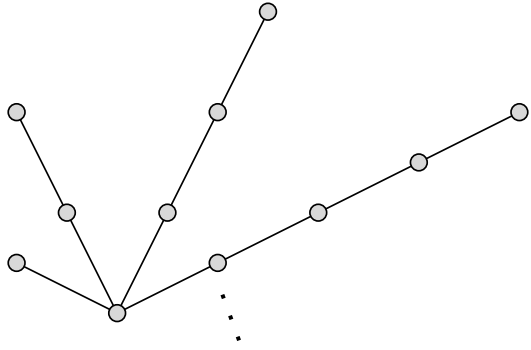, width=3in}
\end{center}
\caption{The tree $T_{\omega}$.}\label{fig:T_omega}
\end{figure}

To further explore the capture-time ordinals, we consider the following example. Define the tree $T_{\omega} = (V(T_\omega), E(T_\omega))$ by setting $$V(T_\omega) = \{r\} \cup \{x_{i, j} \colon 0<i,
j < \omega, j \le i\} \quad \textrm{and}$$ $$E(T_\omega) = \big \{ \{r, x_{i, 1}\} \colon 0<i < \omega\} \big\} \cup \big\{\{x_{i,j}, x_{i,j+1}\} \colon 0<i,j < \omega, j <i \big\}.$$ In other words,
for each positive integer $n$ attach a path of length $n $ to the root vertex {r}. The tree $T_\omega$ is suggested in Figure~\ref{fig:T_omega}.
The tree $T_\omega$ possesses no ray and so is cop-win. However, the
robber may choose as starting position any given end-vertex, making the
capture-time unbounded. We now make this precise in the setting of
CR-ordinals. Note that
$$\eta(x_{i(j+k)},x_{ij}) = i-j  \enskip \textrm{for any~} 0<i,j,k <
\omega \textrm{~such that~} j+k \le i.$$ Hence, we have that
$$\eta(x_{ij}, r) =  i \enskip \textrm{for any~} 0<i,j < \omega, j \le i.$$

For any $v \in V(T_\omega)$ there exists an $i < \omega$ such that $v
\le_i r$. Further, for any $i < \omega$ there exists a $v \in
V(T_\omega)$ such that $v \le_i r$ does not hold. Thus, we have that
$\eta(r) = \omega$.

It is straightforward to derive that for any $0<i < \omega$, $\eta(r, x_{i1}) = \omega$. Further, for any vertex $v \in V(T_\omega)$, $v \le_\omega x_{i1}$. Actually, regardless of the robber's
moves, the cop can move to $r$. Hence, for any $0<i < \omega$ we have that $\eta(x_{i1}) = \omega$.

In addition, for $0<i, j < \omega$, $j \le i$, $\eta(r, x_{ij}) = \omega + j - 1$. Note that for any vertex $v \in V(T_\omega)$ we have $v \le_{\omega + j - 1} x_{ij}$. In addition, for some $v$,
$\eta(v, x_{ij})$ may be smaller than $\omega + j - 1$. Hence, $\eta(x_{ij}) = \omega + j - 1$.
\smallskip

By the above observations, we find that for all positive integers $t$, there are vertices $x$ with $\eta(x) =\omega + t$. Hence, we have the perhaps surprising conclusion that $\rho(T_\omega) =
\omega + \omega = \omega \cdot 2$. Note that $\eta(T_{\omega}) = \omega$, while $\theta(T_\omega) = \{r\} \cup \{x_{i1} \colon 0<i < \omega\}$.

An analogous argument works for any infinite tree obtained by attaching any cardinal number of finite paths to a common root, provided the path lengths are unbounded. Such trees may be uncountable.
This gives examples of uncountable graphs with $\eta(G)= \omega\cdot 2$ (see also the construction in the next section).

\section{The classification of CR-ordinals for trees}
The main result
of the paper is the following theorem.  Define $$\Lambda_T =
\big\{\rho(T) \colon T \textrm{~is a cop-win tree}\big\}.$$
\begin{theorem}{\label{thm:tree}} A CR-ordinal for a cop-win tree is either a finite ordinal, or of the form $\alpha + \omega$, where $\alpha$ is a limit ordinal. In particular,
 $$\Lambda_T = \omega \cup \big\{\alpha+\omega \colon \alpha \textrm{~is a limit ordinal}\big\}.$$
\end{theorem}

The proof of Theorem~\ref{thm:tree} is divided into two parts: necessity (see Lemmas~\ref{lem:finite} and \ref{lem:infinite}) and sufficiency (see Lemma~\ref{lem:suf}). For necessity, the main idea
of the proof is to determine the relationship between the values of $\eta(T)$, $\rho(T)$ and the radius of the given tree. For sufficiency, we give a transfinite construction of trees whose
CR-ordinals have values in $\Lambda_T$.

We first present a lemma which is essentially a part of folklore. We
omit the proof as it is straightforward.

\begin{lemma} \label{lem:finite}
Let $T$ be a tree with finite radius. Then $\eta(T)$ equals the
radius of $T$, and $\rho(T)$ equals the diameter of $T$. In
particular, $\rho(T)$ is a finite ordinal.
\end{lemma}

The next lemma complements the above for the necessity part of
Theorem~\ref{thm:tree}.

\begin{lemma} \label{lem:infinite}
If $T$ is a cop-win tree with infinite radius, then $\eta(T)$ is
infinite and $\rho(T) = \eta(T) + \omega$.
\end{lemma}
In the proof we heavily use that in a tree, an optimal cop strategy is
always to go directly towards the robber using the unique geodesic
connecting the two players' vertices. Observe that this greedy strategy of the cops may not
be optimal for other graph classes. See the results on the game of
Zombies and Survivors \cite{zomb,fhmp}.

\begin{proof}
As $\eta(u, v) \ge d(u, v)$ for any two vertices $u, v \in V(T)$,
we derive immediately that $\eta(T)$ is infinite.

We claim that the set $\theta(T)$ induces a subtree of $T$ with diameter at most two.
 Now, if $\theta(T)$ contains more than one vertex, let $u$ and $v$ be distinct vertices in $\theta(T)$.  If $u$
and $v$ are not adjacent, then let $w$ be any vertex on the unique $uv$-path. Note that vertex $w$ has the property that for any vertex $x \in V(T)$ at least one of $ux$-path or $vx$-path contains
vertex $w$. Therefore, vertex $w$ is also contained in the set $\theta(T)$.  Hence, the subgraph induced by $\theta(T)$ in $T$ is connected.

By rechoosing $u$ and $v$ if necessary, suppose for a contradiction
that the length of the $uv$-path is three; say the path is $uu'v'v$.
Let $T_u$, $T_v$ be the two disjoint subtrees obtained from $T$ by
removing edge $\{u', v'\}$, such that $u \in V(T_u)$ and $v \in
V(T_v)$, respectively. Notice that $u' \in \theta(T)$ implies that
either $\eta(v', u') = \eta(T)$, or for each ordinal $\alpha <
\eta(T)$, there exists a vertex $w \in V(T_u)$ such that $\eta(w, u')
> \alpha$. The second case, however, implies that $\eta(u', v') \ge
\eta(T)$. As $v' \in \theta(T)$ as well, we have that $\eta(v', u') =
\eta(T)$ or $\eta(u', v') = \eta(T)$. But then we have that $$\eta(u)
\ge \eta(v', u) > \eta(T) \mbox{ or } \eta(v) \ge \eta(u', v) >
\eta(T),$$ which contradicts that both $u$ and $v$ are contained in
the set $\theta(T)$. Hence, the claim concerning $\theta(T)$ follows.

By induction we have that for any vertex $v \in V(T) \setminus
\theta(T)$, we have that
$$\eta(v) = \eta(T) + d(v, \theta(T)).$$
As the radius of the tree is infinite while the diameter of the
subgraph induced by $\theta(T)$ is finite, for any $n < \omega$ there
exists some vertex $v \in V(T)$ such that $d(v, \theta(T)) > n$.
Hence, we derive that $\rho(T) = \eta(T) + \omega$.
\end{proof}

By elementary ordinal arithmetic, we may assume that $\eta$ in the sum
$\eta + \omega$ is a limit ordinal. This follows since every successor
ordinal is of the form $\alpha + k$, where $\alpha$ is a limit
ordinal, and $k$ is a finite ordinal.

\medskip

To prove Theorem~\ref{thm:tree}, it is enough to construct a family of trees with CR-ordinals taking all values in the set $\omega \cup \big\{\alpha+\omega \colon \alpha \textrm{~is a limit
ordinal}\big\}$. Finding examples of trees with CR-ordinals equalling all the finite, non-zero ordinals is straightforward. For this, consider the family of finite paths $\{P_n\}_{n\ge 2}$:
Lemma~\ref{lem:finite} implies that $\eta(P_n) = \lceil \frac{n-1}{2} \rceil$ and $\rho(P_n) = n-1.$

We now turn to our construction in the infinite case. We construct a
family $\{T_\alpha \colon \alpha \in \ON \}$, such
that for any $\alpha \in \ON$ we have that $\eta(T_\alpha)
= \alpha$. This construction, in light of Lemma~\ref{lem:infinite},
will complete the proof of Theorem~\ref{thm:tree}.  \smallskip

The construction is based on the operation of \emph{summing rooted
trees}. The basic idea is to form a new root, then append trees to the
root by new edges. To be precise, suppose that $\{(T_i,r_i): i \in
\alpha\}$ is a set of disjoint rooted trees indexed by the ordinal
$\alpha$. Form the rooted tree $\bigoplus_{i \in \alpha} (T_i,r_i)$ by
adding a new vertex $r$ that is joined to each of the $r_i$.

We construct our examples by transfinite induction. Let $S_1 = (K_1,r_1)$ with $r_1$ equalling the single vertex. For any ordinal $\alpha > 1$, assume that all the rooted trees
$(S_{\alpha},r_{\alpha})$ are defined. Let $S_{\alpha}$ be the rooted tree $\bigoplus_{i < \alpha} (S_i,r_i)$, whose root we denote $r_{\alpha}$.

\begin{figure}[htbp]
\begin{center}
\epsfig{clip, figure=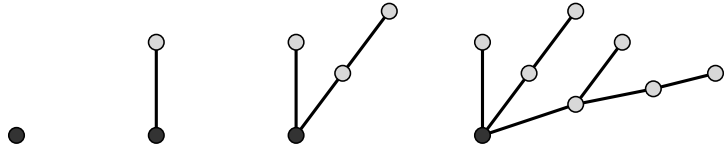, width=4in}
\end{center}
\caption{The rooted trees $S_1$, $S_2$, $S_3$ and $S_4$, with roots coloured black.}\label{fig:RootedTrees}
\end{figure}

See the Figure~\ref{fig:RootedTrees} for the first four trees in the
family $\{ (S_{\alpha},r_{\alpha}): \alpha \in \ON \}$. For
simplicity, we refer to these as $S_{\alpha}$.

As an aside, by Lemma~\ref{lem:finite}, we have that for $0<n < \omega$, $\eta(S_{n+1}) = n$ and $\rho(S_{n+1}) = 2n-1$.

\begin{lemma}\label{lem:suf}
For $\alpha \in \ON$ we have that $\eta(S_{\alpha+1}) = \alpha$.
\end{lemma}
Note that taking $T_\alpha = S_{\alpha+1}$ for $\alpha \in \ON
\setminus \omega$ constructs the desired family.

\begin{proof}
 We prove a slightly stronger statement;
namely that for any $\alpha \in \ON$ there exists a~vertex $v \in
V(S_{\alpha+1})$ such that $$\eta(S_{\alpha+1}) = \eta(r_{\alpha+1}) =
\eta(v, r) = \alpha$$ and, when $\alpha$ is a limit ordinal, there is
$\eta(S_{\alpha}) = \eta(r_{\alpha}) = \alpha$.

We proceed by transfinite induction on $\alpha$, with case $\alpha = 1$ being trivially true. Assume that for some $\alpha \in \ON$ and for all $\beta < \alpha$ one has $\eta (S_{\beta+1}) = \beta$,
and $\eta(r_{\beta+1}) = \eta(v_{\beta+1}, r_{\beta+1}) =\beta$, for some $v_{\beta+1} \in V(S_{\beta+1})$ and for all limit ordinals $\beta < \alpha$ there is $\eta(S_{\beta}) = \eta(r_{\beta}) =
\beta$.

Let $\alpha = \beta + 1$ be a successor ordinal. By construction, we see that $(S_{\alpha+1},r_{\alpha+1})$ consists of the two copies of $(S_{\beta+1},r_{\beta+1})$ on disjoint vertex sets, say $(S,
r)$ and $(S', r')$ joined by the edge $\{r,r'\}$, with $r_{\alpha+1} = r$.

Consider the capture relation in $S_{\alpha+1}$.
Note that every non-leaf vertex of a tree is a cut vertex; hence,
the cop can forbid the robber to enter any component (that is, subtree) except
the one occupied by the robber. Therefore, using the induction
hypothesis, we have that $\eta(u, r) \le \beta$ for all $u \in V(S)$
and that there exists $v \in V(S')$ such that $\eta(v, r') = \beta$,
while $\eta(u, r') \le \beta$ for all $u \in V(S')$. Since $r'$ is a
neighbour of $r$ we have that $\eta(v, r) \le \beta+1$ for any $v \in
V(S_{\alpha+1})$, while $\eta(r', r) = \beta+1$.

Observe that
$\eta(r', v) \ge \eta(r', r)$ for any $v \in V(S)$ and, by the
symmetry, for any $v \in V(S')$ we have that $\eta(r, v) \ge \eta(r,
r') = \beta + 1$. Hence, $$\eta(S_{\alpha+1}) = \eta(r) = \eta(r', r)
= \beta+1=\alpha.$$

Now let $\alpha$ be a limit ordinal. Using the induction hypothesis it is straightforward to see that in $(S_{\alpha}, r_\alpha) = \bigoplus_{\beta < \alpha} (S_\beta,r_\beta)$, for any ordinal
$\beta < \alpha$ there exists a neighbour of $r_\alpha$, namely $r_{\beta+1}$, such that $\eta(r_{\beta+1}, r_\alpha) = \beta+1$. On the other hand, for any $v \in V(S_{\alpha}) \setminus
\{r_\alpha\}$ we have that $v \in V(S_\beta)$, for some $\beta <\alpha$; hence, there is $\eta(v, r) \le \beta+1 < \alpha$. Therefore, we have that $\eta(r_\alpha) = \alpha$, while for any vertex $v
\in V(S_\alpha)$ we have that $\eta(v, r_\alpha) < \alpha$.  Moreover, we find that $\eta(r_\alpha, v) \ge \alpha$, for $v \in V(S_{\alpha}) \setminus \{r_\alpha\}$. This follows since a robber on
$r_\alpha$ may choose to escape to any neighbour of $r_\alpha$ except at most the one leading towards vertex $v$, on which the cop resides. Hence, $\eta(S_\alpha) = \alpha = \eta(r_\alpha)$.

To finish this case consider tree $(S_{\alpha+1}, r_{\alpha+1})$ as two copies of $(S_\alpha,r_\alpha)$ on disjoint vertex sets, say $(S, r)$ and $(S', r')$ joined by the edge $\{r,r'\}$, with
$r_{\alpha+1} = r$. By considering the capture relation on $S_{\alpha+1}$, we see that $\eta(r', r) = \alpha$, while for any vertex $v \in V(S_{\alpha+1})$ there is $\eta(v, r) \le \alpha$. Of course
there is $\eta(r', v) \ge \eta(r', r)$ for any $v \in V(S)$ and, by the symmetry, for any $v \in V(S')$ we have that $\eta(r, v) \ge \eta(r, r') = \alpha$.

Hence, we have that $$\eta(S_{\alpha+1}) = \eta(r) = \eta(r', r) = \alpha,$$ and the proof follows.
\end{proof}

\section{General cop-win graphs}

Define $$\Lambda = \big\{\rho(G) \colon G \textrm{~is a cop-win
graph}\big\}.$$ Of course, $\Lambda_T \subseteq \Lambda$. However, the
containment is strict as we describe in the following example.

\begin{figure}[htbp]
\begin{center}
\epsfig{clip, figure=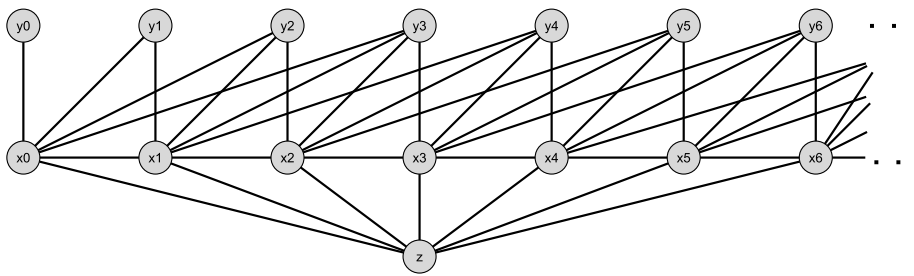, width=5in}
\end{center}
\caption{The Polat graph.}\label{fig:G1}
\end{figure}

Consider the following cop-win graph introduced in~\cite{polat}, which we refer to as the \emph{Polat graph.} Let $X = \{x_n \colon n < \omega\}$, $Y = \{y_n \colon n < \omega\}$, and $Z = \{z\}$ be
disjoint sets of vertices. Let $G$ be the graph defined by $V(G) = X \cup Y \cup Z$ and $$E(G) = \bigcup_{n < \omega} \Big(\{x_n, x_{n+1}\}, \{x_n, z\}, \{x_n, y_{n}\}, \{x_n, y_{n+1}\}, \{x_n,
y_{n+2}\}, \{x_n, y_{n+3}\}\Big).$$ The Polat graph is suggested in Figure~\ref{fig:G1}.
By direct checking, we have that $\rho(G) = \omega+1$, $\eta(G) = \omega$ and $\theta(G) = X \cup Z$ (we omit the details here). Hence, the Polat graph witnesses the fact that $\omega +1 \in \Lambda \setminus \Lambda_T$.

We may modify the construction of the Polat graph in several ways to obtain other CR-ordinals. For example, we may add a finite path to the vertex $z$ obtaining graphs with CR-ordinal equalling
$\omega+i$, for $0<i < \omega$. By forming a sum (akin to the rooted tree sum described in the previous section) of these generalized Polat graphs, for $j < \omega, j >1$, there exist cop-win graphs
with CR-ordinal $\omega \cdot j + (i+j)$.
\section{Problems}
The main open problem we consider is to classify which ordinals belong to
$\Lambda$. Define
$$
\Upsilon = \big\{\omega \cdot i + (i + j) \colon i,j < \omega\} \cup \big\{\alpha+\omega \colon \alpha \textrm{~is a limit ordinal}\big\}.$$ Our family of cop-win graphs derived from the Polat graph
supports the assertion that $\Lambda = \Upsilon$. We leave this as an open problem. Some ordinals do not seem possible to attain as a CR-ordinal. For example, is $\omega \in \Lambda$? Note that this
question is answered negatively if $\Lambda = \Upsilon$. Observe that $\rho(G)$ is well-defined for any (not necessarily cop-win) graph. Another question, therefore, is to classify the ordinals
$\rho(G)$.

We mention in closing that the paper~\cite{nw} does not distinguish between cardinals and ordinals (though its results are correct when the distinction is made). Further, we do not know the best
upper bound on the tower of relations $\{\le_\alpha\}_{\alpha < \aleph_{\gamma+1}}$ for a graph $G$ of cardinality $\aleph_{\gamma}$.

\end{document}